\documentclass[9pt,technote]{IEEEtran}
\usepackage{fancyhdr}

\usepackage{amsfonts}
\usepackage{mathrsfs}
\usepackage{colortbl}

\usepackage{bm}
\usepackage{amsmath,amssymb,amsthm,graphicx}
\usepackage{subfigure,epstopdf}  %%%these two packages used together

\newtheorem{Assumption}{\bf Assumption}
\newtheorem{Lemma}{\bf Lemma}
\newtheorem{Theorem}{\bf Theorem}

\newtheorem{Remark}{\bf Remark}

\newcommand{\RR}{\mathbb{R}}

\newcommand{\EQ}{\begin{eqnarray}}

\newcommand{\EN}{\end{eqnarray}}

\newcommand{\EQQ}{\begin{eqnarray*}}

\newcommand{\ENN}{\end{eqnarray*}}

\newenvironment{Proof}{\noindent{\em Proof:\/}}{\hfill $\Box$\par}

\begin{document}

\title{Adaptive Leader-Following Consensus for Uncertain Euler-Lagrange Systems under Directed Switching Networks}

\author{Tao~Liu~and~Jie~Huang
\thanks{This work has been supported by the Research Grants Council of the
Hong Kong Special Administration Region under grant No. 412612.}
\thanks{Tao Liu and Jie Huang  are with the Department of Mechanical and Automation
Engineering, The Chinese University of Hong Kong Shatin, N.T., Hong
Kong. E-mail: tliu2@mae.cuhk.edu.hk,jhuang@mae.cuhk.edu.hk}
}

\maketitle

\begin{abstract}
The leader-following consensus problem for multiple Euler-Lagrange systems was studied recently by the adaptive
distributed observer approach under the assumptions that the leader system is neurally stable and the communication network
is jointly connected and undirected.
In this paper, we will study the same problem without assuming that the leader system is neutrally stable, and the communication network is
undirected. The effectiveness of this new result will be illustrated by an example.
\end{abstract}

%\begin{IEEEkeywords}
%IEEEtran, journal, \LaTeX, paper, template.
%\end{IEEEkeywords}

% \fi
\IEEEpeerreviewmaketitle

\section{Introduction}
Consensus, as a fundamental problem of cooperative control, has received significant attention over the past decade
\cite{JadbaLin03}, \cite{OlfatiFax07}, \cite{RenBeard05}, \cite{Tuna09}.
There are two types of consensus problems, i.e., leaderless consensus and leader-following consensus.
The leaderless consensus problem aims to design a distributed control law to make the states/outputs of all agents synchronize to each other,
while the leader-following consensus problem attempts to drive the states/outputs of all agents to a prescribed trajectory generated by a leader system.
Euler-Lagrange (EL) systems is an important class of nonlinear systems, that  models a large class of mechanical systems including robotic manipulators
and rigid bodies \cite{LewisAbdal93}, \cite{SlotineLi91}. The consensus problem for multiple EL systems has  been extensively investigated.

The leader-following consensus problem for multiple EL systems was first considered in \cite{ChungSlotine09} assuming that all followers have access to the leader.
The same problem was further studied in \cite{MeiRen11} under the assumption that the communication network of the multiple EL systems is static, undirected and connected, and
in \cite{ChenLewis11}, \cite{NunoOrtega11} under the assumption that the communication network of the multiple EL systems is static and connected.

More recently, the leader-following consensus problem for multiple EL systems subject to jointly connected switching communication network was studied
\cite{CaiHuang14}, \cite{CaiHuang16}. Specifically, by employing a distributed observer, a distributed adaptive state feedback control law was synthesized to
solve the  leader-following consensus problem for multiple EL systems under a set of standard assumptions in \cite{CaiHuang14}.
A drawback of the distributed observer in \cite{CaiHuang14} is that the system matrix of the leader has to be used by all followers, which may not be realistic in some applications.
This drawback was overcome in \cite{CaiHuang16} by replacing the distributed observer with a so-called adaptive distributed observer, which is capable of providing the estimated system
matrix of the leader to  all followers. Thus the control law in \cite{CaiHuang16} does not require  the system matrix of the leader  be used by all the followers.
Nevertheless, the success of \cite{CaiHuang16} was obtained at two other costs. First, it required that the leader system be neurally stable, which precludes the frequently used ramp signal. Second, it assumed that the  communication network was undirected, which also limited the scope of the applications of the result in
\cite{CaiHuang16}.

In this paper, we will offer two improvements over the main result in \cite{CaiHuang16}. That is, we will obtain the same result as in
\cite{CaiHuang16} using the adaptive distributed observer approach but without assuming that the leader system is neutrally stable and the  communication network is undirected.
For this purpose, we need to first strengthen the result on the adaptive distributed observer \cite{CaiHuang16} so that it applies to unbounded leader's signal in polynomial form. Then we will establish our main result using this strengthened version of the adaptive distributed observer.

In what follows, we will adopt the following notation. $\bm{1}_{N}$ denotes an $N$ dimensional column vector whose components are all $1$.
$\otimes$ denotes the Kronecker product of matrices. $\|x\|$ denotes the Euclidean norm of a vector $x$ and $\|A\|$ denotes the induced norm of a matrix $A$
by the Euclidean norm. $\lambda_{\text{max}}(A)$ and $\lambda_{\text{min}}(A)$ denote the maximum and the minimum eigenvalues of a matrix $A$, respectively.
For $X_{i} \in \RR^{n_{i} \times p}, i=1, \ldots, m$, col$\left(X_{1}, \ldots, X_{m}\right)=\left[X_{1}^{T}, \ldots, X_{m}^{T}\right]^{T}$.
We call a time function $\sigma:[0, +\infty) \mapsto \mathcal{P}=\{1,2,\ldots,n_{0}\}$ a piecewise constant switching signal if there exists a sequence
$\{t_{i}, i=0,1,2,\ldots \}$ satisfying $t_{0}=0, t_{i+1}-t_{i}\ge \tau_{0}$ for some positive constant $\tau_{0}$,
such that, for all $t\in [t_{i},t_{i+1})$, $\sigma(t)=p$ for some $p \in \mathcal{P}$.
$n_{0}$ is some positive integer.
$\mathcal{P}$ is called the switching index set;
$t_{i}$ is called the switching instant and $\tau_{0}$ is called the dwell time.

\section{Problem Formulation and Assumptions}
Consider $N$ EL systems described by the following dynamic equations:
\begin{equation}\label{eq1}
  M_{i}(q_{i})\ddot{q}_{i}+C_{i}(q_{i},\dot{q}_{i})\dot{q}_{i}+G_{i}(q_{i})=\tau_{i}, \quad i=1, \ldots, N
\end{equation}
where $q_{i}, \dot{q}_{i} \in \RR^{n}$ are the generalized position and velocity vectors, respectively; $M_{i}(q_{i}) \in \RR^{n \times n}$ is the positive definite inertia matrix;
$C_{i}(q_{i},\dot{q}_{i})\dot{q}_{i} \in \RR^{n} $ is the Coriolis and centripetal forces vector; $G_{i}(q_{i}) \in \RR^{n} $ is the gravity vector,
and $\tau_{i} \in \RR^{n}$ is the generalized forces vector.

It is well known that  the EL systems have the following two properties:

\textbf{\emph{Property 1:}} $\dot{M}_{i}(q_{i})-2C_{i}(q_{i},\dot{q}_{i})$ is skew symmetric.

\textbf{\emph{Property 2:}} For all $x,y \in \RR^{n}$,
\[M_{i}(q_{i})x+C_{i}(q_{i},\dot{q}_{i})y+G_{i}(q_{i})=Y_{i}(q_{i},\dot{q}_{i},x,y)\Theta_{i}\]
where $Y_{i}(q_{i},\dot{q}_{i},x,y) \in \RR^{n \times p}$ is a known regression matrix
and $\Theta_{i} \in \RR^{p}$ is a constant vector consisting of the uncertain parameters of (\ref{eq1}).

Like in \cite{CaiHuang14},~\cite{CaiHuang16}, let $q_{0} \in \RR^{n}$ denote the desired generalized position vector, which is assumed to be generated by the following exosystem:
\begin{equation}\label{eq2}
  \dot{v}=Sv, \qquad q_{0}=Cv
\end{equation}
where $v \in \RR^{m}$ and $S \in \RR^{m \times m}, C \in \RR^{n \times m}$ are constant matrices. Without loss of generality, we assume
the pair $(C, S)$ is observable.

We view the system composed of (\ref{eq1}) and (\ref{eq2}) as a multi-agent system of $(N+1)$ agents
with (\ref{eq2}) as the leader and $N$ subsystems of (\ref{eq1}) as followers.
Given systems (\ref{eq1}), (\ref{eq2}) and a piecewise constant switching signal $\sigma(t)$,
we can define a switching digraph $\bar{\mathcal{G}}_{\sigma(t)}=(\bar{\mathcal{V}},\bar{\mathcal{E}}_{\sigma(t)})$\footnote{See Appendix for a summary on digraph.}
with $\bar{\mathcal{V}}=\{0,1,\ldots, N\}$ and $\bar{\mathcal{E}}_{\sigma(t)} \subseteq \bar{\mathcal{V}} \times \bar{\mathcal{V}}$ for all $t\ge 0$.
Here, node $0$ is associated with the leader system (\ref{eq2}) and node $i, \ i=1,\ldots, N$, is associated with the $i$th subsystem of (\ref{eq1}).
For $i=0,1,\ldots,N, j=1,\ldots,N$, $(i,j) \in \bar{\mathcal{E}}_{\sigma(t)}$ if and only if $\tau_{j}$
can use the state of agent $i$ for control at time instant $t$.
As a result, our control law has to satisfy the communication constraint described by the digraph $\bar{\mathcal{G}}_{\sigma(t)}$.
Such a control law is called a distributed control law.

Our problem is described as follows.

\textbf{\emph{Problem Description:}} Given systems (\ref{eq1}), (\ref{eq2}) and a switching digraph $\bar{\mathcal{G}}_{\sigma(t)}$,
find a distributed state feedback control law of the following form:
\begin{align}\label{eq3}
  \tau_{i} & = f_{i}\left(q_{i},\dot{q}_{i},\varphi_{i},\varphi_{j}-\varphi_{i},j \in \bar{\mathcal{N}}_{i}(t) \right)    \notag  \\
  \dot{\varphi}_{i} & = g_{i} \left (\varphi_{i},\varphi_{j}-\varphi_{i},j \in \bar{\mathcal{N}}_{i}(t) \right), \quad i=1,\ldots,N
\end{align}
where $\bar{\mathcal{N}}_{i}(t)$ denotes the neighbor set of agent $i$ at time $t$,
such that, for $i=1, \ldots, N$, and for any initial conditions $v(0)$, $q_{i}(0)$ and $\dot{q}_{i}(0)$,
$q_{i}(t)$ and $\dot{q}_{i}(t)$ exist for all $t \ge 0$ and satisfy
\begin{equation}\label{eq4}
  \lim_{t \to +\infty}(q_{i}(t) - q_{0}(t))=0,  \quad  \lim_{t \to +\infty}(\dot{q}_{i}(t)-\dot{q}_{0}(t))=0.
\end{equation}

Some assumptions for the solvability of the above problem are listed below.

\begin{Assumption}\label{ass01}
None of the eigenvalues of $S$ have positive real parts.
\end{Assumption}

\begin{Assumption}\label{ass04}
$\dot{q}_{0}$ is bounded.
\end{Assumption}

\begin{Assumption}\label{ass02}
There exist positive constants $k_{\underline{m}}$, $k_{\overline{m}}$, $k_{c}$, $k_{g}$,
such that, for $i=1,\ldots,N$, $k_{\underline{m}}I_{n} \le M_{i}(q_{i}) \le  k_{\overline{m}} I_{n}$,
$\|C_{i}(q_{i},\dot{q}_{i})  \|  \le k_{c}\|\dot{q}_{i}\|$, and $\| G_{i}(q_{i}) \|  \le k_{g}$.
\end{Assumption}

\begin{Assumption}\label{ass03}
There exists a subsequence $\{i_{k}\}$, $k=0,1,2,\ldots$, of $\{i : i=0,1,2,\ldots\}$ with $t_{i_{k+1}}-t_{i_{k}} < \epsilon$
for some positive $\epsilon$ such that every node $i$, $i=1,\ldots,N$, is reachable from node $0$ in the union digraph
$\bigcup_{j=i_{k}}^{i_{k+1}-1} \bar{\mathcal{G}}_{\sigma(t_{j})}$.
\end{Assumption}

\begin{Remark}\label{re01.1}
Assumption \ref{ass01} allows the  generalized position vector $q_{0}$ of the leader system (\ref{eq2}) to be a polynomial in $t$ and thus is much more general
than the assumption that the leader system is neutrally stable required in~\cite{CaiHuang16}.
Assumption \ref{ass04} is more restrictive than Assumption \ref{ass01}. However, it still allows
the  generalized position vector $q_{0}$ of the leader system (\ref{eq2})
to be a ramp function, which is not allowed in ~\cite{CaiHuang16}.
\end{Remark}

\begin{Remark}\label{re01}
Assumption \ref{ass03} is called the jointly connected condition~\cite{JadbaLin03} and is perhaps the mildest condition on a switching network
since it allows the network to be disconnected at any time instant.
\end{Remark}

\section{Main Results}
Let us first recall the adaptive distributed observer introduced in \cite{CaiHuang16}.
For this purpose, let $\bar{\mathcal{A}}_{\sigma(t)}=[a_{ij}(t)]_{i,j=0}^{N}$ denote the weighted adjacency matrix of $\bar{\mathcal{G}}_{\sigma(t)}$.
Then, for each agent of (\ref{eq1}), we define a dynamic compensator as follows:
\begin{align}\label{eq5}
  \dot{S}_{i} &= \mu_{1} \sum_{j=0}^{N}a_{ij}(t)(S_{j}-S_{i})    \notag   \\
  \dot{\eta}_{i} & = S_{i}\eta_{i} + \mu_{2} \sum_{j=0}^{N}a_{ij}(t) (\eta_{j} - \eta_{i}), \quad i=1,\ldots,N
\end{align}
where $S_{i} \in \RR^{m \times m}, S_{0}=S, \eta_{i} \in \RR^{m}, \eta_{0}=v$, $\mu_{1}$ and $\mu_{2}$ are any positive constants.

Furthermore, let $\mathcal{G}_{\sigma(t)}=(\mathcal{V},\mathcal{E}_{\sigma(t)})$ denote the subgraph of $\bar{\mathcal{G}}_{\sigma(t)}$,
where $\mathcal{V}=\{1,\ldots, N\}$ and $\mathcal{E}_{\sigma(t)} \subseteq \mathcal{V} \times \mathcal{V}$ is obtained from $\bar{\mathcal{E}}_{\sigma(t)}$
by removing all the edges between node $0$ and the nodes in $\mathcal{V}$. Let $\mathcal{L}_{\sigma(t)}$ be the Laplacian of $\mathcal{G}_{\sigma(t)}$.
Then, putting $\eta=\text{col}(\eta_{1},\ldots,\eta_{N})$, $\hat{\eta}=\eta-\bm{1}_{N} \otimes v$, $\hat{S}_{i}=S_{i}-S$, $\hat{S}=\text{col}(\hat{S}_{1},\ldots,\hat{S}_{N})$
and $\hat{S}_{d}=\text{block diag} \left\{\hat{S}_{1},\ldots,\hat{S}_{N} \right\}$, we can write (\ref{eq5}) into the following compact form:
\begin{align}\label{eq6}
  \dot{\hat{S}} & = - \mu_{1}\left(H_{\sigma(t)} \otimes I_{m} \right)\hat{S}      \notag \\
  \dot{\hat{\eta}} & = \left (I_{N} \otimes S -\mu_{2} (H_{\sigma(t)} \otimes I_{m}) \right )\hat{\eta} + \hat{S}_{d} \eta
\end{align}
where $H_{\sigma(t)}=\mathcal{L}_{\sigma(t)}+\text{diag}\left\{a_{10}(t),\ldots,a_{N0}(t)\right\}$.

Now, let us establish the following result.
\begin{Lemma}\label{le01}
Under Assumptions \ref{ass01} and \ref{ass03}, for any $\mu_{1},\mu_{2}>0$, and for any initial conditions $\hat{S}(0)$ and $\hat{\eta}(0)$,
we have
\begin{equation}\label{}
  \lim_{t\to +\infty} \hat{S}(t)=0
\end{equation}
exponentially, and
\begin{equation}\label{eq7}
  \lim_{t\to +\infty} \hat{\eta}(t)=0
\end{equation}
asymptotically.
\end{Lemma}

\begin{Proof}
By Corollary 4 of~\cite{SuHuang12A}, for any $\mu_{1}>0$, the origin of the $\hat{S}$-subsystem of (\ref{eq6}) is exponentially stable.
That is to say, $\lim_{t\to +\infty} \hat{S}(t)=0$, exponentially. Thus, we only need to prove (\ref{eq7}).
Denote $A(t)=\left (I_{N} \otimes S -\mu_{2} (H_{\sigma(t)} \otimes I_{m}) \right )$ and $F(t)= \hat{S}_{d}(t)(\bm{1}_{N} \otimes v)$.
Then, the second equation of (\ref{eq6}) is equivalent to
\begin{equation}\label{eq7y}
  \dot{\hat{\eta}}  = A(t)\hat{\eta} + \hat{S}_{d}(t)\hat{\eta} + F(t).
\end{equation}

Since  $\hat{S}_{d}(t)$ converges to zero exponentially, there exist $\alpha_{1}>0$ and $\lambda_{1}>0$ such that
\begin{equation}\label{}
  \| \hat{S}_{d}(t)  \|  \le  \alpha_{1}\| \hat{S}_{d}(0)  \| e^{-\lambda_{1}t}.
\end{equation}
Note that
\begin{equation}\label{}
  \| (\bm{1}_{N} \otimes v) \| \le \| (I_{N} \otimes e^{St}) \| \,  \| (\bm{1}_{N} \otimes v(0))   \|.
\end{equation}
Under Assumption \ref{ass01}, there exists a polynomial $p(t)$ such that
\begin{equation}\label{}
   \| (I_{N} \otimes e^{St}) \| \le p(t).
\end{equation}
Then,
\begin{align}\label{}
  \|F(t)\| & \le  \| \hat{S}_{d}(t)  \| \, \| (\bm{1}_{N} \otimes v) \|         \notag   \\
           & \le \alpha_{1} \| \hat{S}_{d}(0)  \| \, \| (\bm{1}_{N} \otimes v(0)) \| p(t) e^{-\lambda_{1}t}    \notag \\
           & \le \alpha_{2} \| \hat{S}_{d}(0)  \| \, \| (\bm{1}_{N} \otimes v(0)) \| e^{-\lambda_{2}t}
\end{align}
for some $\alpha_{2} >0$ and $\lambda_{1} > \lambda_{2} >0$. Thus, $F(t)$ also converges to zero exponentially.

By Lemma 2 of~\cite{SuHuang12A}, under Assumptions \ref{ass01} and \ref{ass03}, for any $\mu_{2}>0$,
the origin of the linear switched system
\begin{equation}\label{eq7x}
  \dot{\hat{\eta}}  = A(t)\hat{\eta}
\end{equation}
is exponentially stable.
Let $\Phi(\tau,t)\hat{\eta}$ be the solution of (\ref{eq7x}) that starts at $(t,\hat{\eta})$.
Define
\begin{equation}\label{}
  P(t) = \int_{t}^{+\infty} \Phi(\tau,t)^{T}Q\Phi(\tau,t)d \tau
\end{equation}
where $Q$ is some constant positive definite matrix. Clearly, $P(t)$ is continuous for all $t \ge 0$.
Since the equilibrium point $\hat{\eta}=0$ of (\ref{eq7x}) is exponentially stable, we have
\begin{equation}\label{}
  \| \Phi(\tau,t)\|  \le \alpha_{3} e^{-\lambda_{3}(\tau-t)}, \quad \forall \tau \ge t \ge 0
\end{equation}
for some $\alpha_{3}>0$ and $\lambda_{3}>0$.
It can be easily verified that $c_{1}\|\hat{\eta}\|^{2} \le \hat{\eta}^{T}P(t) \hat{\eta} \le  c_{2}\|\hat{\eta}\|^{2}$
for some positive constants $c_{1}$ and $c_{2}$. Hence $P(t)$ is positive definite and bounded.
Thus, we can assume that $\|P(t)\| \le c_{3}$ for any $t \ge 0$ with $c_{3}$ being some positive constant.

On the other hand, since $A(t)$ is continuous on intervals $[t_{i},t_{i+1}), i=0,1,2,\ldots$, we have,
for $t \in [t_{i}, t_{i+1})$,  $i=0,1,2,\ldots$,
\begin{equation}\label{}
  \frac{\partial}{\partial t} \Phi(\tau,t)=-\Phi(\tau,t) A(t), \quad \Phi(t,t)=I_{m}.
\end{equation}
Then we have
\begin{align}\label{}
  \dot{P}(t) & = \int_{t}^{+\infty}  \Phi(\tau,t)^{T}Q   \left(\frac{\partial}{\partial t}\Phi(\tau,t)  \right )d \tau    \notag \\
             &  \quad +   \int_{t}^{+\infty}   \left(\frac{\partial}{\partial t}\Phi(\tau,t)^{T} \right ) Q   \Phi(\tau,t)  d \tau  -Q     \notag  \\
   &= - \int_{t}^{+\infty}  \Phi(\tau,t)^{T}  Q  \Phi(\tau,t)  d \tau  A(t)      \notag  \\
   &   \quad - A(t)^{T}\int_{t}^{+\infty}  \Phi(\tau,t)^{T}  Q  \Phi(\tau,t)  d \tau -Q   \notag  \\
   &= -P(t)A(t)-A(t)^{T}P(t)-Q.
\end{align}

Let $U(t)=\hat{\eta}^{T} (t) P(t)\hat{\eta} (t)$. Then,
along the trajectory of (\ref{eq7y}), for any $t \in [t_{i},t_{i+1})$ with $i=0,1,2,\ldots$, we have
\begin{align}\label{eq7z}
  \dot{U}(t)& = \hat{\eta}^{T} \left( \dot{P}(t) + A(t)^{T}P(t) + P(t)A(t)  \right )\hat{\eta}    \notag  \\
                                       & \quad   + 2 \hat{\eta}^{T}P(t)\hat{S}_{d}(t)\hat{\eta} + 2 \hat{\eta}^{T}P(t)F(t)            \notag  \\
  & = - \hat{\eta}^{T}Q\hat{\eta} + 2 \hat{\eta}^{T}P(t)\hat{S}_{d}(t)\hat{\eta} + 2 \hat{\eta}^{T}P(t)F(t)    \notag  \\
  & \le - \hat{\eta}^{T}Q\hat{\eta} + 2c_{3} \| \hat{S}_{d}(t) \|  \,  \| \hat{\eta} \|^{2} +  2 \hat{\eta}^{T}P(t)F(t)  \notag  \\
  &  \le  -\lambda_{\text{min}}(Q)\| \hat{\eta} \|^{2} + 2c_{3} \| \hat{S}_{d}(t) \| \, \| \hat{\eta} \|^{2}     \notag \\
  &  \quad  + \frac{\|P(t)\|^{2}}{\varepsilon} \| \hat{\eta} \|^{2} + \varepsilon \|F(t)\|^{2}   \notag  \\
  &  \le  - \left( \lambda_{\text{min}}(Q) -  2c_{3} \| \hat{S}_{d}(t) \|  -  \frac{c_{3}^{2}}{\varepsilon}    \right) \| \hat{\eta} \|^{2}    \notag   \\
  & \quad   + \varepsilon \|F(t)\|^{2}.
\end{align}
Choose $\varepsilon=\frac{2c_{3}^{2}}{\lambda_{\text{min}}(Q)}$.
Then, since $\hat{S}_{d}(t)$ converges to zero exponentially, there exists some positive integer $l$, such that
\begin{equation}\label{}
  \left( \lambda_{\text{min}}(Q) -  2c_{3} \| \hat{S}_{d}(t) \|  -  \frac{c_{3}^{2}}{\varepsilon}    \right) >0, \quad \forall t \ge t_{l}.
\end{equation}
Thus, we have
\begin{equation}\label{}
  \dot{U}(t) \le \varepsilon \|F(t)\|^{2},  \quad \forall t \ge t_{l}
\end{equation}
which implies
\begin{equation}\label{}
  {U}(t) \le   {U}(t_l) + \varepsilon \int_{t_l}^t \|F(\tau)\|^{2} d \tau,  \quad \forall t \ge t_{l}.
\end{equation}
Since $F(t)$  converges to zero exponentially, $\lim_{t \to +\infty}U(t)$ exists and is finite.
Thus,  we conclude that $U(t)$ is bounded over $t \geq 0$ and hence the solution
$\hat{\eta}(t)$ of (\ref{eq7y}) is also bounded over $t \geq 0$.

In addition, for any $t \in [t_{i}, t_{i+1}), i=0,1,2,\ldots$, we have
$\ddot{U}(t)$ is bounded over $[0,+\infty)$ since $\hat{\eta}$, $\dot{\hat{\eta}}$, $P(t)$, $\dot{P}(t)$,  $\hat{S}_{d}(t)$, $\dot{\hat{S}}_{d}(t)$
$F(t)$ and $\dot{F}(t)$, are all bounded over $[0,+\infty)$.

Thus, $U (t)$ satisfies the three conditions of Lemma $1$ of~\cite{SuHuang12B}. As a result, $\dot{U}(t) \to 0$ as $t \to +\infty$, which in turn
implies that the solution $\hat{\eta}(t)$ of (\ref{eq7y}) converges to zero asymptotically.
Hence the proof is completed.
\end{Proof}

\begin{Remark}\label{re01x}
Since $\dot{U} $ is only piecewise continuous over $[0,+\infty)$, instead of using Barbala's lemma, we have to use
Lemma $1$ of~\cite{SuHuang12B} to conclude $\dot{U} (t) \to 0$ as $t \to +\infty$.
\end{Remark}

\begin{Remark}\label{re03}
As a result of Lemma \ref{le01}, under Assumptions \ref{ass01} and \ref{ass03},
for any $\mu_{1}, \mu_{2} >0$, and $i=1,\ldots,N$,
\EQ
  \lim_{t \to +\infty} \left( S_{i}(t)-S    \right) =0   \label{eq22a}
\EN
\EQ
  \lim_{t \to +\infty} \left( \eta_{i}(t)-v(t)   \right) =0. \label{eq22b}
\EN
That is why (\ref{eq5}) is called the adaptive distributed observer of the leader system (\ref{eq2}).
Moreover, let $\eta_{di}=\mu_{2}\sum_{j=0}^{N}a_{ij}(t) (\eta_{j} - \eta_{i})$. Then,  (\ref{eq22b}) implies
\EQ \label{eq22c}
 \lim_{t \rightarrow \infty} \eta_{di} (t) = 0.
\EN
Since
\begin{align*}
 \dot{\eta}_{i}-\dot{v} &=S_{i}\eta_{i}+\eta_{di}-Sv   \\
                         &=S_{i}(\eta_{i}-v)+\hat{S}_{i}v+\eta_{di}, \\
\end{align*}
we have
\EQ \label{eq22d}
 \lim_{t \rightarrow \infty} (\dot{\eta}_{i}-\dot{v} ) = 0, \quad i=1,\ldots,N.
 \EN
\end{Remark}

\begin{Remark}\label{re02}
The adaptive distributed observer for the leader system (\ref{eq2}) was first developed in Lemma 2 of \cite{CaiHuang16}
under the assumptions that all the eigenvalues of the matrix $S$ are semi-simple with zero real parts and
the digraph ${\cal G}_{\sigma (t)}$ is undirected.
Lemma $2$ of \cite{CaiHuang16} was strengthened recently by Lemma $4.1$ of \cite{LiuHuang16},
which  removed  the assumption that the digraph ${\cal G}_{\sigma (t)}$ is undirected.
Here, Lemma \ref{le01} further replaced the neutral stability assumption on the matrix $S$ required in \cite{CaiHuang16} and \cite{LiuHuang16} with Assumption \ref{ass01}.
As a result, we can handle signals in polynomial form.
\end{Remark}

%\begin{Remark}\label{re04}
% and .
%Then,
%\begin{equation}\label{eq8}
%  \lim_{t \to +\infty}(\xi_{i}(t)-q_{0}(t))=\lim_{t \to +\infty}C(\eta_{i}(t)-v(t))=0
%\end{equation}
%and
%\begin{align}\label{eq9}
%\lim_{t \to +\infty}(\dot{\xi}_{i}(t)- \dot{q}_{0}(t)) &= \lim_{t \to +\infty} C\left(  S_{i}(t)\eta_{i}(t) +  \eta_{di}(t) -Sv(t)    \right)   \notag  \\
%& = \lim_{t \to +\infty} C ((S_{i}(t)-S)\eta_{i}(t) + \eta_{di}(t)   \notag  \\
%&  \qquad \qquad  +S(\eta_{i}(t)-v(t)) )   \notag \\
%& =0.
%\end{align}
%\end{Remark}

Next, like in \cite{CaiHuang16}, we will synthesize an adaptive distributed control law utilizing the adaptive distributed observer as follows.

Let $\xi_{i}=C\eta_{i}$ and
\begin{equation}\label{eq10}
  \dot{q}_{ri}= CS_{i}\eta_{i} -\alpha(q_{i}-\xi_{i})
\end{equation}
where $\alpha$ is a positive constant. Then,
\begin{equation}\label{eq11}
  \ddot{q}_{ri}=C\left( \dot{S}_{i}\eta_{i} + S_{i} \dot{\eta}_{i}\right)  -\alpha(\dot{q}_{i}-\dot{\xi}_{i})   .
\end{equation}
By Property $2$, there exists a known matrix $Y_{i}=Y_{i}(q_{i},\dot{q}_{i},\ddot{q}_{ri},\dot{q}_{ri})$
and an unknown constant vector $\Theta_{i}$ such that
\begin{equation}\label{eq12}
  Y_{i}\Theta_{i}= M_{i}(q_{i})\ddot{q}_{ri}+C_{i}(q_{i},\dot{q}_{i})\dot{q}_{ri}+G_{i}(q_{i}).
\end{equation}

Let
\begin{equation}\label{eq13}
  s_{i}=\dot{q}_{i}-\dot{q}_{ri}.
\end{equation}
Then, we define our control law as follows:
\begin{align}
  \tau_{i} & =-K_{i}s_{i} + Y_{i} \hat{\Theta}_{i}            \label{eq14}  \\
  \dot{\hat{\Theta}}_{i} & = - \Lambda_{i}^{-1}Y_{i}^{T}s_{i}    \label{eq15}   \\
  \dot{S}_{i} &= \mu_{1} \sum_{j=0}^{N}a_{ij}(t)(S_{j}-S_{i})       \\
  \dot{\eta}_{i} & = S_{i}\eta_{i} + \mu_{2} \sum_{j=0}^{N}a_{ij}(t) (\eta_{j} - \eta_{i}), \quad i=1,\ldots,N      \label{eq17}
\end{align}
where $\hat{\Theta}_{i} \in \RR^{p}$, $K_{i}$ and $\Lambda_{i}$ are positive definite matrices.

Now, we are ready to present our main result.

\begin{Theorem}\label{thrm01}
Given systems (\ref{eq1}), (\ref{eq2}) and a switching digraph $\bar{\mathcal{G}}_{\sigma(t)}$,
under Assumptions \ref{ass04} to \ref{ass03},
the problem is solvable by a distributed state feedback control law composed of (\ref{eq14})-(\ref{eq17}).
\end{Theorem}

%\subsection{Bounded $q_{0}$, $\dot{q}_{0}$ and $\ddot{q}_{0}$.}
%In this case, the generalized position, speed and acceleration vectors are \mbox{all} bounded, and we can solve the problem without Assumption \ref{ass02}.

\begin{Proof}
First note that, under Assumption \ref{ass04}, the leader system also satisfies Assumption \ref{ass01}.
Next, from (\ref{eq10}) and (\ref{eq13}), we have
\begin{equation}\label{eq25}
  \dot{q}_{i}+\alpha (q_{i}-\xi_{i})=s_{i}+CS_{i}\eta_{i}
\end{equation}
where $CS_{i}\eta_{i}=C(\dot{\eta_{i}} - \eta_{di}) = \dot{\xi}_{i}-C\eta_{di}$.
Subtracting $\dot{\xi}_{i}$ on both sides of (\ref{eq25}) gives
\begin{equation}\label{eq25.01}
  (\dot{q}_{i}-\dot{\xi}_{i})+\alpha (q_{i} - \xi_{i}) = u_i
\end{equation}
where $u_i = s_{i} -C\eta_{di}$.  Since $\alpha >0$,   (\ref{eq25.01}) is a stable first order linear system
in $(q_{i} - \xi_{i})$ with  input $u_i$. If $u_i$ decays to zero as $t$ tends to infinity, then
both $(q_{i} - \xi_{i})$ and $(\dot{q}_{i} - \dot{\xi}_{i})$ decay to zero as $t$ tends to infinity. As a result,
by (\ref{eq22b}), (\ref{eq22d}) and the following identities
\begin{align}\label{}
  q_{i}(t) - q_{0}(t) & =(q_{i}(t)-\xi_{i}(t)) + C (\eta_{i}(t)- v (t))   \notag \\
  \dot{q}_{i}(t) - \dot{q}_{0}(t) & =(\dot{q}_{i}(t)-\dot{\xi}_{i}(t)) + C (\dot{\eta}_{i}(t)- \dot{v}(t))
\end{align}
the proof is completed.

By (\ref{eq22c}), under Assumptions \ref{ass04} and \ref{ass03}, $\eta_{di} (t) \to 0$ as $t \to +\infty$. We only need to show $s_{i} (t) \to 0$ as $t \to +\infty$.
To this end, substituting (\ref{eq14}) into (\ref{eq1}) gives
\begin{equation}\label{eq18}
  M_{i}(q_{i})\ddot{q}_{i}+C_{i}(q_{i},\dot{q}_{i})\dot{q}_{i}+G_{i}(q_{i})=-K_{i}s_{i} + Y_{i} \hat{\Theta}_{i}
\end{equation}
and subtracting $Y_{i}\Theta_{i}$ on both sides of (\ref{eq18}) gives
\begin{align}\label{}
  M_{i}(q_{i})\ddot{q}_{i}+C_{i}(q_{i},\dot{q}_{i})\dot{q}_{i}-M_{i}(q_{i})\ddot{q}_{ri} - C_{i}(q_{i},\dot{q}_{i})\dot{q}_{ri} \notag \\
  = -K_{i}s_{i} + Y_{i} \tilde{\Theta}_{i}
\end{align}
where $\tilde{\Theta}_{i}=\hat{\Theta}_{i}-\Theta_{i} $. Then, by (\ref{eq13}), we have
\begin{equation}\label{eq20}
    M_{i}(q_{i})\dot{s}_{i}+C_{i}(q_{i},\dot{q}_{i})s_{i} = -K_{i}s_{i} + Y_{i}\tilde{\Theta}_{i}.
\end{equation}

Let $x=\text{col}(x_{1},\ldots,x_{N})$ for $x=q, \dot{q}, s, \dot{s}, \tilde{\Theta}$,
and $X= \text{block diag}\{ X_{1},\ldots,X_{N}\}$ for $X=K, Y, \Lambda^{-1}$.
Then (\ref{eq20}) and (\ref{eq15}) can be written as
\begin{align}
  M(q)\dot{s} &=-C(q,\dot{q})s-Ks+Y\tilde{\Theta}    \label{eq21}   \\
  \dot{\tilde{\Theta}} &= -\Lambda^{-1}Y^{T}s        \label{eq22}
\end{align}
where
\[M(q) = \text{block diag} \left \{ M_{1}(q_{1}),\ldots,M_{N}(q_{N}) \right \}  \]
\[C(q,\dot{q})=\text{block diag} \left \{C_{1}(q_{1},\dot{q}_{1}),\ldots,C_{N}(q_{N},\dot{q}_{N})\right \}.\]

Define
\begin{equation}\label{eq23}
  V=\frac{1}{2}\left( s^{T}M(q)s + \tilde{\Theta}^{T}\Lambda \tilde{\Theta}    \right).
\end{equation}
By (\ref{eq11}) and (\ref{eq13}), $s(t)$ is differentiable on each interval $[t_{i}, t_{i+1})$,
$i=0,1,2,\ldots$, so is $\dot{V}(t)$. Noticing that $\dot{M}_{i}(q_{i})-2C_{i}(q_{i},\dot{q}_{i})$ is skew symmetric gives
\begin{align}\label{eq24}
  \dot{V} & =s^{T}M(q)\dot{s}+\frac{1}{2}s^{T}\dot{M}(q)s+\tilde{\Theta}^{T}\Lambda\dot{\tilde{\Theta}}       \notag   \\
   & = s^{T}\left (-C(q,\dot{q})s-Ks+Y\tilde{\Theta} \right)+\frac{1}{2}s^{T}\dot{M}(q)s+\tilde{\Theta}^{T}\Lambda\dot{\tilde{\Theta}}   \notag \\
   & = -s^{T}Ks+s^{T}Y\tilde{\Theta}-\tilde{\Theta}^{T}\Lambda\Lambda^{-1}Y^{T}s    \notag \\
   & = -s^{T}Ks \le 0.
\end{align}

Since $V (t)$ and $\dot{V} (t)$ are piecewise continuous over $[0,+\infty)$, we cannot use Barbala's lemma to conclude $\dot{V} (t) \to 0$ as $t \to +\infty$.
We need to use Corollary $1$ of~\cite{SuHuang12B} to conclude $\lim_{t \to +\infty} \dot{V}(t) =0$, which implies
$\lim_{t \to +\infty} s(t)=0$. For this purpose, we need to show that
 there exists a positive number $\gamma $ such that
\begin{equation}\label{eq26}
  \sup_{t_{i} \le t \le t_{i+1}, \ i=0,1,2,\ldots}  |\ddot{V}(t)| \le \gamma.
\end{equation}

Since $\ddot{V}(t)=-2s^{T}K\dot{s}$, it suffices to show that both $s$ and $\dot{s}$ are bounded.

Now note that $V(t)$ is continuous,  and $M(q)$ and $\Lambda$ are positive definite,
 (\ref{eq24}) implies that $s$ and $\tilde{\Theta}$ are bounded. Thus, the input $u_i$ in (\ref{eq25.01}) is bounded.

>From (\ref{eq21}), to show $\dot{s}$ is bounded, we need to show $C (q,\dot{q})$ and $Y\tilde{\Theta}$   are  bounded.

We first note that  (\ref{eq25.01}) implies both both $(q_{i} - \xi_{i})$ and $(\dot{q}_{i} - \dot{\xi}_{i})$ are bounded since $u_i$  is bounded.
By (\ref{eq22d}), $\dot{\xi}_{i} =  C \dot{\eta}_{i}$ is bounded since $\dot{q}_{0} = C \dot{v}$ is bounded.
Thus $\dot{q}_{i}$ is bounded, which implies  $C_{i}(q_{i},\dot{q}_{i})$ is bounded under  Assumption \ref{ass02}.

>From (\ref{eq12}), $Y\tilde{\Theta}$ is bounded if both $\dot{q}_{ri}$ and $\ddot{q}_{ri}$ are bounded.
Since we have already shown that $s_i$ and $\dot{q}_{i}$ are bounded, we have $\dot{q}_{ri}$ is bounded by (\ref{eq13}).

We now show $\ddot{q}_{ri}$ is bounded using (\ref{eq11}).  In fact,
\begin{align}\label{eq25.3}
C\dot{S}_{i}\eta_{i}&=C\dot{\hat{S}}_{i}\eta_{i}=C\dot{\hat{S}}_{i}v+C\dot{\hat{S}}_{i}(\eta_{i}-v) \\
  CS_{i}\dot{\eta}_{i}&=C\hat{S}_{i}\dot{\eta}_{i}+CS\dot{\eta}_{i}    \notag \\
                      &=C\hat{S}_{i}\dot{v}+C\hat{S}_{i}(\dot{\eta}_{i}-\dot{v})
                      +CS\dot{v}+CS(\dot{\eta}_{i}-\dot{v})  \notag\\
                      &=C\hat{S}_{i}\dot{v}+C\hat{S}_{i}(\dot{\eta}_{i}-\dot{v})
                      +\ddot{q}_{0}+CS(\dot{\eta}_{i}-\dot{v}).
\end{align}
Thus, $\ddot{q}_{ri}$ is bounded since, by Remark \ref{re03}, under Assumptions \ref{ass04} and \ref{ass03}, every term on the right hand side of (\ref{eq11}) is bounded.
Thus, (\ref{eq26}) is satisfied.
The proof is completed by invoking Corollary $1$ of~\cite{SuHuang12B}.
\end{Proof}
%$\hat{S}_{i}$, $\dot{\hat{S}}_{i}$ decay to zero as $t$ tends to infinity, and $v$ is a polynomial function, which may contain sinusoidal components.
%Therefore, both $\dot{q}_{ri}$ and $\ddot{q}_{ri}$ are bounded.
%Then, by (\ref{eq12}), $Y_{i}$ is bounded.

%Since $M(q)$ is positive definite, $M(q)^{-1}$ exists and is bounded.
%From (\ref{eq21}), $\dot{s}$ is also bounded.
%
%
%%At last, substituting (\ref{eq10}) into (\ref{eq13}) gives
%%\begin{equation}\label{eq27}
%%  \dot{q}_{i}+\alpha (q_{i} - \xi_{i}) =s_{i} + CS_{i}\eta_{i}.
%%\end{equation}
%%Subtracting $\dot{\xi}_{i}$ on both sides of (\ref{eq27}) gives
%%\begin{equation}\label{eq28}
%%  (\dot{q}_{i}-\dot{\xi}_{i})+\alpha (q_{i} - \xi_{i}) =s_{i} -C\eta_{di}.
%%\end{equation}
%Now, we can view (\ref{eq25.01}) as a stable first order linear system in $(q_{i} - \xi_{i})$
%with a bounded input, and this input tends to zero as $t \to +\infty$.
%Thus, both $(q_{i} - \xi_{i})$ and $(\dot{q}_{i}-\dot{\xi}_{i})$ are bounded over $t \ge 0$
%and will decay to zero.

\begin{Remark}\label{05}
If we strengthen Assumption \ref{ass01} to the one that the leader system is neutrally stable as assumed in~\cite{CaiHuang16},
then the generalized position vector $q_{0}$ as well as its derivative of any degree is bounded.
In this case, Assumption \ref{ass04} is satisfied automatically.
Furthermore, $\xi_{i}$, $\dot{\xi}_{i}$ are bounded from (\ref{eq22b}) and (\ref{eq22d}),
which implies that $q_{i}$, $\dot{q}_{i}$ are bounded. Thus, Assumption \ref{ass02} is also satisfied automatically.
It is worth mentioning that even in this case, we have extended the result of~\cite{CaiHuang16} from undirected communication networks to directed communication networks.
\end{Remark}

\section{An Example}

In this section, we consider a group of four EL systems, each of which describes a two-link manipulator whose
motion equation is taken from~\cite{LewisAbdal93}:
\begin{equation*}\label{}
  M_{i}(q_{i})\ddot{q}_{i}+C_{i}(q_{i},\dot{q}_{i})\dot{q}_{i}+G_{i}(q_{i})=\tau_{i}, \quad i=1, 2, 3, 4,
\end{equation*}
where $q_{i}=\text{col}(\theta_{i1},\theta_{i2})$ and
\begin{align*}\label{}
  M_{i}(q_{i}) & = \left(
                     \begin{array}{cc}
                       a_{i1}+a_{i2}+2a_{i3}\cos \theta_{i2} &  a_{i2}+a_{i3} \cos \theta_{i2} \\
                       a_{i2}+a_{i3} \cos \theta_{i2} &  a_{i2} \\
                     \end{array}
                   \right)           \\
  C_{i}(q_{i},\dot{q}_{i}) & = \left(
                                 \begin{array}{cc}
                                   -a_{i3}(\sin \theta_{i2})\dot{\theta}_{i2} &    -a_{i3}(\sin \theta_{i2})(\dot{\theta}_{i1}+\dot{\theta}_{i2}) \\
                                   a_{i3}(\sin \theta_{i2})\dot{\theta}_{i1} & 0 \\
                                 \end{array}
                               \right)   \\
 G_{i}(q_{i}) & = \left(
                 \begin{array}{c}
                   a_{i4}g\cos \theta_{i1} + a_{i5}g\cos (\theta_{i1}+\theta_{i2}) \\
                   a_{i5}g\cos (\theta_{i1}+\theta_{i2})  \\
                 \end{array}
               \right)
\end{align*}
with $\Theta_{i}=\text{col}(a_{i1}, a_{i2}, a_{i3}, a_{i4}, a_{i5})$. Then, Assumption \ref{ass02} is satisfied.

Let the leader's signal be as follows:
\[
q_{0}(t)=\left(
           \begin{array}{c}
             1+t+\cos t + \sin t \\
             1+t+\cos t -\sin t \\
           \end{array}
         \right).
\]
Then this leader's signal can be produced by the following leader system:
\begin{align*}
  \dot{v} &= Sv=\left(
         \begin{array}{cccc}
           0 & 1 & 0 & 0 \\
           0 & 0 & 0 & 0 \\
           0 & 0 & 0 & 1 \\
           0 & 0 & -1 & 0 \\
         \end{array}
       \right) v \\
  q_0 &= Cv=\left(
          \begin{array}{cccc}
            1 & 0 & 1 & 0 \\
            1 & 0 & 0 & 1 \\
          \end{array}
        \right) v
\end{align*}
%\EQQ
%\dot{v} &=& Sv=\left(
%         \begin{array}{cccc}
%           0 & 1 & 0 & 0 \\
%           0 & 0 & 0 & 0 \\
%           0 & 0 & 0 & 1 \\
%           0 & 0 & -1 & 0 \\
%         \end{array}
%       \right) v \\
% q_0 &=& C v = \left(
%          \begin{array}{cccc}
%            1 & 0 & 1 & 0 \\
%            1 & 0 & 0 & 1 \\
%          \end{array}
%        \right) v
%\ENN
with initial condition $v(0)=\bm{1}_{4}$. It can be verified that the pair $(C, S)$ is observable and  Assumption \ref{ass04} is satisfied.

\begin{figure}%[!htb]
  \centering
%  \subfigure[$\bar{\mathcal{G}}_{1}$]{
%    %\label{g1} %% label for first subfigure
%    \includegraphics[width=1.3in]{figure/g1.eps}}
%  \hspace{0.5in}
  \subfigure[$\bar{\mathcal{G}}_{1}$]{
   % \label{g2} %% label for second subfigure
    \includegraphics[width=0.9in]{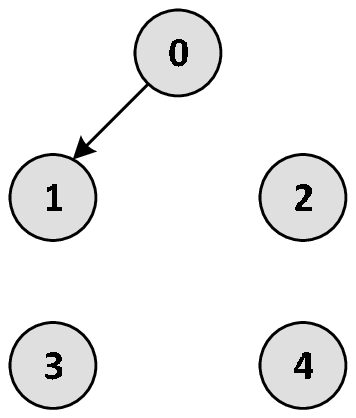}}
      \hspace{0.4in}
  \subfigure[$\bar{\mathcal{G}}_{2}$]{
   % \label{g3} %% label for second subfigure
    \includegraphics[width=0.9in]{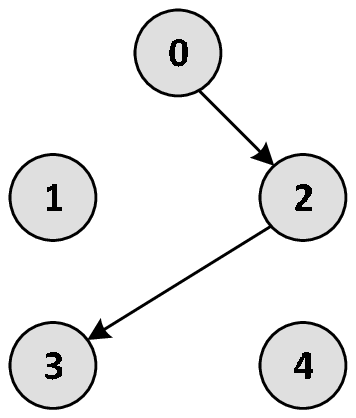}}
          \hspace{0.4in}
  \subfigure[$\bar{\mathcal{G}}_{3}$]{
   % \label{g2} %% label for second subfigure
    \includegraphics[width=0.9in]{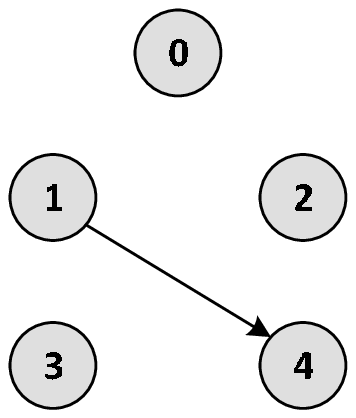}}
     \hspace{0.4in}
  \subfigure[$\bar{\mathcal{G}}_{4}$]{
   % \label{g3} %% label for second subfigure
    \includegraphics[width=0.9in]{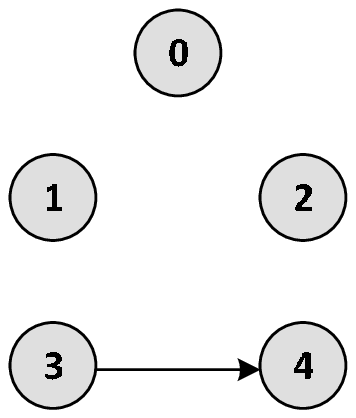}}
  \caption{Switching topology $\bar{\mathcal{G}}_{\sigma(t)}$ with
$\mathcal{P}=\{1,2,3,4\}$}
  \label{g} %% label for entire figure
\end{figure}

Let the switching digraph $\bar{\mathcal{G}}_{\sigma(t)}$ be dictated by the following switching signal:
\begin{equation}\label{}
   \sigma(t)=
\begin{cases}
1, & \textrm{if}\quad sT_{0} \le t <(s+ \frac{1}{4})T_{0} \\
2, & \textrm{if}\quad (s+ \frac{1}{4})T_{0} \le t < (s+ \frac{1}{2})T_{0} \\
3, & \textrm{if}\quad (s+ \frac{1}{2})T_{0} \le t < (s+ \frac{3}{4})T_{0}  \\
4, & \textrm{if}\quad (s+ \frac{3}{4})T_{0} \le t < (s+ 1)T_{0}
\end{cases}
\end{equation}
where $T_{0}=2$, and $s=0,1,2,\ldots$.
The four digraphs $\bar{\mathcal{G}}_{i}, i=1,2,3,4$, are described by Figure \ref{g}
where node $0$ is associated with the leader and the other \mbox{nodes} are associated with the followers.
It can be seen that Assumption \ref{ass03} is satisfied even though $\bar{\mathcal{G}}_{\sigma(t)}$ is disconnected at any time $t \ge 0$.

%\emph{Case A: Bounded $q_{0}$.}
%
%
%Let
%\[  C = \left(
%          \begin{array}{cccc}
%            0 & 1 & 1 & 0 \\
%            0 & 1 & 0 & 1 \\
%          \end{array}
%        \right).
% \]
%Then
%\[
%q_{0}(t)=\left(
%           \begin{array}{c}
%             1+\cos t + \sin t \\
%             1+\cos t -\sin t \\
%           \end{array}
%         \right).
%\]
%
%
%
%\emph{Case B: Unbounded $q_{0}$.}
%
%
%Let
%\[  C = \left(
%          \begin{array}{cccc}
%            1 & 0 & 1 & 0 \\
%            1 & 0 & 0 & 1 \\
%          \end{array}
%        \right).
% \]
%Then
%\[
%q_{0}(t)=\left(
%           \begin{array}{c}
%             1+t+\cos t + \sin t \\
%             1+t+\cos t -\sin t \\
%           \end{array}
%         \right).
%\]

According to Theorem \ref{thrm01}, we can design a control law in the form described by (\ref{eq14})-(\ref{eq17})
with the following design parameters: $\mu_{1}=\mu_{2}=10$, $\alpha=10$, $K_{i}=20I_{2}$, $\Lambda_{i}=0.2I_{5}$, for $i=1,2,3,4$.
We let $a_{ij}(t)=1$, $i,j=0,1,2,3,4$, whenever $(j,i) \in \bar{\mathcal{E}}_{\sigma(t)}$.
The actual values of $\Theta_{i}$ are given as follows:
\begin{align*}
  \Theta_{1} & =\text{col}(0.64, 1.10, 0.08, 0.64, 0.32) \\
  \Theta_{2} & =\text{col}(0.76, 1.17, 0.14, 0.93, 0.44) \\
  \Theta_{3} & =\text{col}(0.91, 1.26, 0.22, 1.27, 0.58) \\
  \Theta_{4} & =\text{col}(1.10, 1.36, 0.32, 1.67, 0.73).
\end{align*}

\begin{figure}
  \centering
  \includegraphics[width=3.5in,height=4in]{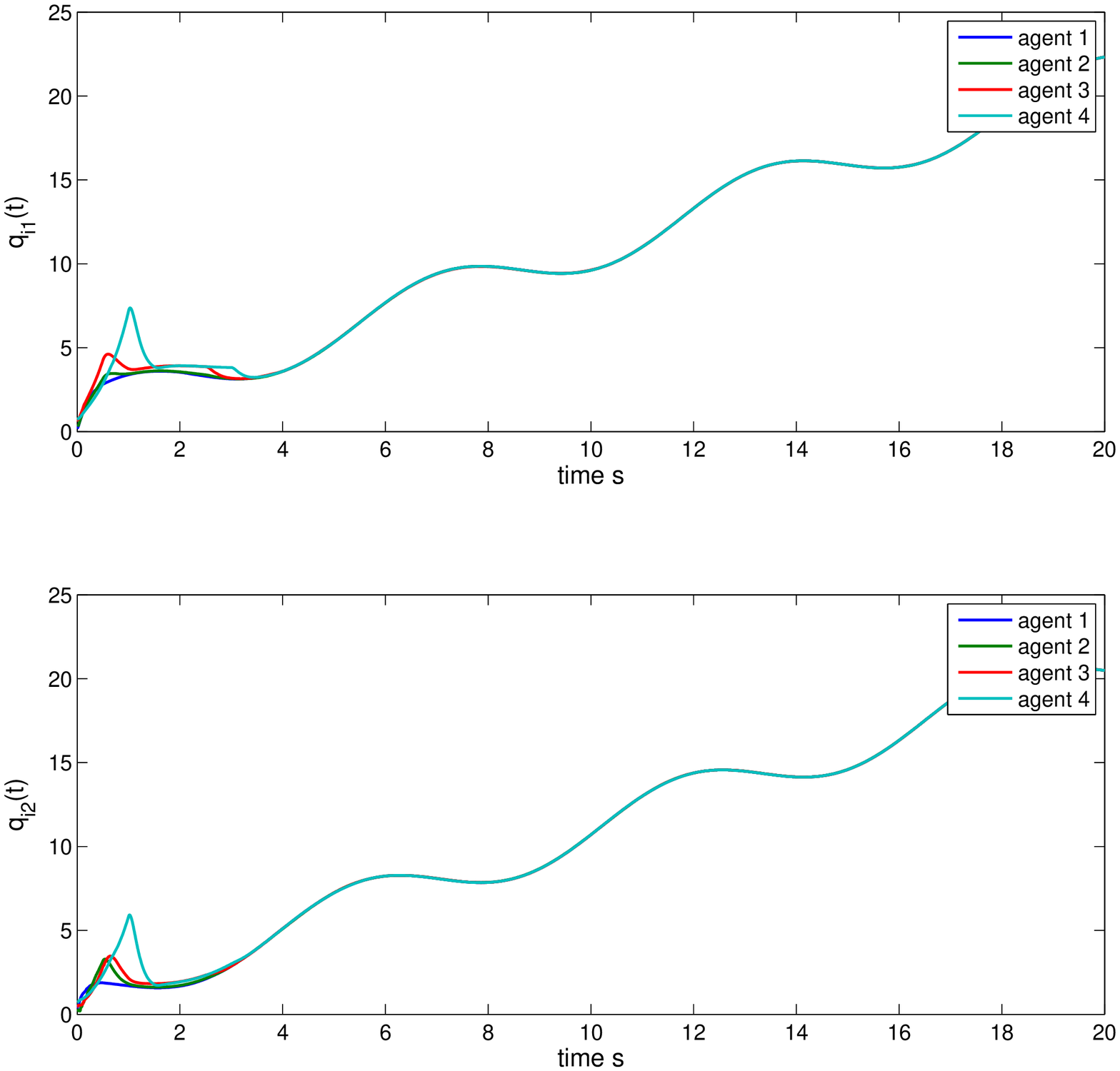}
  \caption{Generalized position of each agent}
\label{position}
\end{figure}

\begin{figure}
  \centering
  \includegraphics[width=3.5in,height=4in]{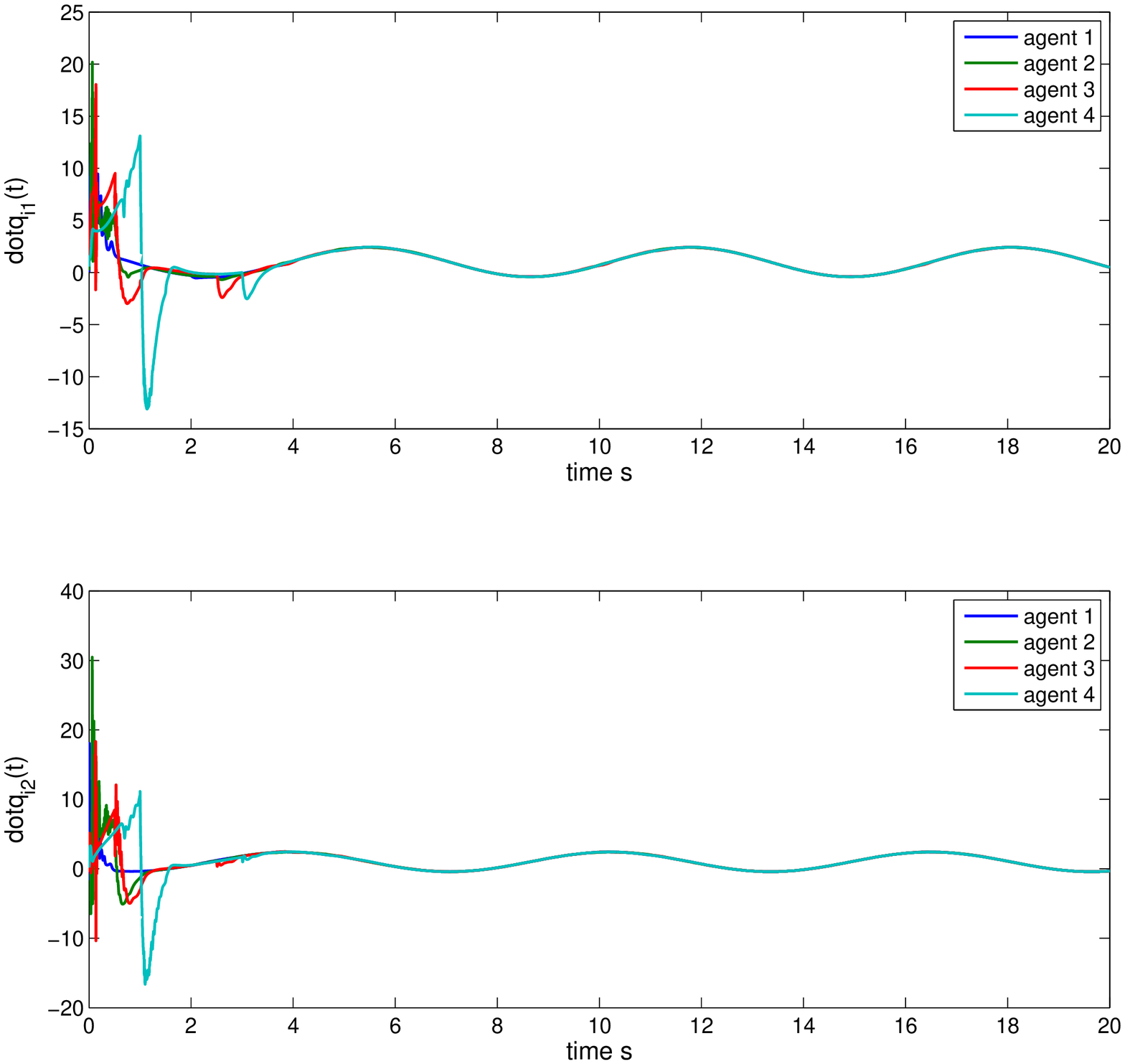}
  \caption{Generalized velocity of each agent}
\label{velocity}
\end{figure}

Simulation is conducted with randomly chosen initial conditions. The trajectories of $q_{i}$ and $\dot{q}_{i}$, $i=1,2,3,4$,
are shown in Figure \ref{position} and Figure \ref{velocity}, respectively.

\section{Conclusion}
In this paper, we have studied the leader-following consensus problem for multiple uncertain Euler-Lagrange systems
under the jointly connected switching network. Due to the employment of the adaptive distributed observer in a strengthened version,
we have removed the assumptions that the leader system is neutrally stable and the communication network is undirected.

\section*{Appendix}
A digraph $\mathcal{G}=(\mathcal{V}, \mathcal{E})$ consists of a finite set of nodes $\mathcal{V}=\{1, \ldots, N\}$
and an edge set $\mathcal{E} \subseteq \mathcal{V} \times \mathcal{V}$. An edge of $\mathcal{E}$ from node $i$ to node $j$
is denoted by $(i,j)$, and node $i$ is called a neighbor of node $j$. Let $\mathcal{N}_{i}=\{j|(j,i) \in \mathcal{E} \}$,
which is called the neighbor set of node $i$.
The edge $(i,j)$ is called undirected if $(i,j) \in \mathcal{E}$ implies $(j,i) \in \mathcal{E}$.
The digraph $\mathcal{G}$ is undirected if every edge in $\mathcal{E}$ is undirected.
If the digraph contains a sequence of edges of the form $(i_{1},i_{2})$, $(i_{2},i_{3})$,
$\ldots$, $(i_{k},i_{k+1})$, then the set $\left \{(i_{1},i_{2}), (i_{2},i_{3}), \ldots, (i_{k},i_{k+1}) \right \}$ is called a directed path of $\mathcal{G}$
from node $i_{1}$ to node $i_{k+1}$ and node $i_{k+1}$ is said to be reachable from node $i_{1}$.
A digraph $\mathcal{G}_{s}=(\mathcal{V}_{s}, \mathcal{E}_{s})$ is called a subgraph of $\mathcal{G}=(\mathcal{V}, \mathcal{E})$
if $\mathcal{V}_{s} \subseteq \mathcal{V} $ and $\mathcal{E}_{s} \subseteq \mathcal{E}\bigcap (\mathcal{V}_{s} \times \mathcal{V}_{s} )$.
Given a set of $n_{0}$ digraphs $\left \{ \mathcal{G}_{i}=(\mathcal{V}, \mathcal{E}_{i}), i=1,\ldots,n_{0} \right \}$, the digraph
$\mathcal{G}=(\mathcal{V}, \mathcal{E})$ where $\mathcal{E}= \bigcup_{i=1}^{n_{0}} \mathcal{E}_{i}$ is called the union of the digraphs $\mathcal{G}_{i}$,
denoted by $\mathcal{G}=\bigcup_{i=1}^{n_{0}} \mathcal{G}_{i}$.

The weighted adjacency matrix of a digraph $\mathcal{G}$ is a nonnegative matrix $\mathcal{A}=[a_{ij}] \in \RR^{N \times N}$, where
$a_{ii}=0$ and $a_{ij} > 0$ if and only if $(j,i)\in \mathcal{E}, i,j=1,\ldots,N$. On the other hand, given a matrix $\mathcal{A}=[a_{ij}] \in \RR^{N \times N}$
satisfying $a_{ii}=0$ and $a_{ij} \ge 0$ for $i \ne j$, we can always define a digraph $\mathcal{G}$ whose weighted adjacency matrix is $\mathcal{A}$.
The Laplacian of $\mathcal{G}$ is then defined as $\mathcal{L}=[l_{ij}] \in \RR^{N \times N}$, where $l_{ii}=\sum_{j=1}^{N}a_{ij}$, $l_{ij}=-a_{ij}$ for $i \ne j$.

Given a piecewise constant switching signal $\sigma: [0, +\infty) \mapsto \mathcal{P}=\{1,2,\dots,n_{0}\} $, and a set of $n_{0}$ digraphs
$\mathcal{G}_{i}=(\mathcal{V}, \mathcal{E}_{i})$, $i=1,\ldots,n_{0}$, with the corresponding weighted adjacency matrices being denoted by $\mathcal{A}_{i}$,
$i=1, \ldots, n_{0}$, we call the time-varying graph $\mathcal{G}_{\sigma(t)}=(\mathcal{V},\mathcal{E}_{\sigma(t)})$ a switching digraph, and denote
the weighted adjacency matrix and the Laplacian of $\mathcal{G}_{\sigma(t)}$ by $\mathcal{A}_{\sigma(t)}$ and $\mathcal{L}_{\sigma(t)}$, respectively.

%\section{Acknowledgements}
%This work has been supported by the Research Grants Council of the Hong Kong Special Administration Region under grant No. 412408,
%and in part by the National Natural Science Foundation of China under grant No. 61174049.

\end{document}